\documentclass[a4paper,12pt]{article}
\usepackage[all]{xy}
\usepackage{amsmath, amsthm}
\usepackage{amssymb,amsfonts}
\usepackage{array}
\usepackage{hyperref}

\title{Vertex and Edge Connectivity of the Zero Divisor graph $\Gamma[\mathbb {Z}_n]$}
\author{B.Surendranath Reddy,Rupali.S.Jain and N.Laxmikanth\\
\textit{surendra.phd@gmail.com},{rupalisjain@gmail.com} and\\
{laxmikanth.nandala@gmail.com} }{}
\date{}
\begin{document}
\theoremstyle{definition}
\newtheorem{definition}{Definition}[section]
\newtheorem{example}[definition]{Example}
\newtheorem{remark}[definition]{Remark}
\newtheorem{observation}[definition]{Observation}
\theoremstyle{plain}
\newtheorem{theorem}[definition]{Theorem}
\newtheorem{lemma}[definition]{Lemma}
\newtheorem{proposition}[definition]{Proposition}
\newtheorem{corollary}[definition]{Corollary}
\newtheorem{AMS}[definition]{AMS}
\newtheorem{keyword}[definition]{keyword}
\maketitle
\section*{Abstract}
The Zero divisor Graph of a commutative ring $R$, denoted by $\Gamma[R]$, is a graph whose vertices are non-zero zero divisors of $R$ and two vertices are adjacent if their product is zero. In this paper we derive the Vertex and Edge Connectivity of the zero divisor graph $\Gamma[\mathbb{Z}_n]$, for any natural number $n$ . We also discuss the minimum degree of the zero divisor graph $\Gamma[\mathbb{Z}_n]$.\\
Keywords:-Zero divisor graph, Vertex connectivity,Edge connectivity,\\minimum degree.
\section{Introduction}
The concept of the Zero divisor graph of a ring $R$ was first introduced by I.Beck[\cite{3}] in 1988 and later on Anderson and Livingston[\cite{2}], Akbari and Mohammadian[\cite{1}]  continued the study of zero divisor graph by considering only the non-zero zero divisors. Mohammad Reza and Reza Jahani[\cite{5}] calculated the energy and Wiener index for the zero divisor graphs  $\Gamma[{\mathbb{Z}_n}]$ for $n=p^2$ and $n=pq$ where $p$ and $q$ are prime numbers. B.Surendranath Reddy, et.al[\cite{6}] considered the zero divisor graph $\Gamma[{\mathbb{Z}_n}]$ for $n=p^3\, \text{and}\, p^2q$, where $p$ and $q$ are prime numbers and derived the standard form of adjacency matrix, spectrum, energy and Wiener index of. The concepts of the Edge and Vertex connectivity  of a graph can be found in [\cite{4}]. In this paper we extend the concepts of  the Edge and Vertex connectivity to the  zero divisor graph $\Gamma[{\mathbb{Z}_n}]$ for any natural number $n$.

 In this article,  section 2, is about the preliminaries and notations related to zero divisor graph of a commutative ring $R$, in section 3, we derive the  Vertex connectivity of a zero divisor graph  $\Gamma[{\mathbb{Z}_n}]$, and in section 4, we calculate the  Edge connectivity of $\Gamma[{\mathbb{Z}_n}]$.In this section,  we also find the  minimum degree  of the zero divisor graph $\Gamma[{\mathbb{Z}_n}]$.
 \section{Preliminaries and Notations}
\begin{definition}{\textbf{Zero divisor Graph}}{\cite{3}}\\
Let R be a commutative ring with unity and $ Z[R]$ be the set of its zero divisors. Then the zero divisor graph of R denoted by  $\Gamma[R]$, is the graph(undirected) with vertex set  $Z^*[R]= Z[R]-\{\mathbf{0}\}$, the non-zero zero divisors of $R$, such that two vertices $v,w \in Z^*[R]$ are adjacent if $vw=0$.
\end{definition}
\begin{definition}{\textbf{Vertex connectivity of a graph}}{\cite{4}}\\
 Let $G$ be a simple graph. The Vertex connectivity of $G$, denoted by $\kappa(G)$, is the smallest number of vertices in $G$ whose deletion from $G$ leaves either a disconnected graph or $K_{1}$.
\end{definition}
\begin{definition}{\textbf{Edge connectivity of a graph}}{\cite{4}}\\
 Let $G$ be a simple graph. The Edge connectivity of $G$, denoted by $\kappa_{e}(G)$, is the smallest number of edges in $G$ whose deletion from $G$ leaves either a disconnected graph or an empty graph.
\end{definition}
\begin{definition}{\textbf{Minimum degree  of a graph }}{\cite{4}}\\
 Let $G$ be a  graph. The Minimum degree of $G$, denoted by $\delta(G)$, is the minimum degree of its vertices.
\end{definition}
\section{Vertex connectivity of of the zero divisor graph $\Gamma[\mathbb{Z}_n]$}
In this section we derive the vertex connectivity of $\Gamma[\mathbb{Z}_n]$ for any natural number $n$. To start with we consider first the zero divisor graph $\Gamma[\mathbb{Z}_n]$ for $n=p^2$.
\begin{theorem}
The vertex connectivity of $\Gamma[\mathbb{Z}_{p^2}]$ is $p-2$.
 \begin{proof}
The vertex set of the zero divisor graph $\Gamma[\mathbb{Z}_{p^2}]$ is 
 $A =\{kp\,|\,k= 1,2,3,....,p-1 \}$ and so $|A|=(p-1)$.\\
 As product of any two vertices is zero, they are adjacent and so the corresponding graph is a complete graph on (p-1) vertices that is, $\Gamma[\mathbb{Z}_{p^2}]=K_{p-1}$.\\
As the graph is complete, deletion of one or less than $p-2$ vertices does not give a disconnected graph. If we delete $p-2$ vertices then only one vertex remains and that gives rise to a disconnected graph.\\
Thus the minimum number of vertices to be deleted is $p-2$.\\
Hence the vertex connectivity of $\Gamma[\mathbb{Z}_{p^2}]$ is $p-2$.\\
Since $\kappa(\Gamma[\mathbb{Z}_{p^2}])= p-2$, the zero divisor graph $\Gamma[\mathbb{Z}_{p^2}]$  is $p-2$ connected.
 \end{proof}
\end{theorem}
Now we consider the case for $n=p^3$  in the following theorem.
\begin{theorem}
The vertex connectivity of $\Gamma[\mathbb{Z}_{p^3}]$ is $p-1$.
 \begin{proof}
consider the zero divisor graph $\Gamma[\mathbb{Z}_{p^3}]$.\\
Here, we divide the elements(vertices) of $\Gamma[\mathbb{Z}_{p^3}]$ into two disjoint sets namely multiples of $p$ but not $p^2$; and the multiples of $p^2$ which are given by
\begin{align*}
 A &=\{kp\,|\,k= 1,2,3,....,p^2-1\, \text{and}\, p\nmid k \}  \\
 B &=\{lp^2\,|\,l= 1,2,3,....,p-1 \} 
\end{align*}
with cardinality  $|A|=p(p-1)$ and $|B|=(p-1)$.\\
Here we note that no two elements of $A$ are adjacent and every element of $A$ is adjacent with every element of $B$. Also every element of $B$  is adjacent with every element of $A$ and $B$. \\
Now if we remove all the vertices of $A$, the graph is still connected because the remaning vertices are from set $B$ which are adjacent with each other.
whereas, if we remove all the vertices from $B$, then the resulting graph with vertices  from  $A$ is disconnected as no two elements of $A$ are adjacent.\\
Also if we leave even one vertex from $B$ and remove remaining, the graph is still connected as every element of $B$ is adjacent with $A$ and $B$.\\
Therefore, the minimum number of vertices to be deleted from $G$ is the number of vertices of $B$.\\
Hence the vertex connectivity $\kappa(\Gamma[\mathbb{Z}_{p^3}])=|B|=p-1$.\\
Since, $\kappa(\Gamma[\mathbb{Z}_{p^3}])=p-1$, the zero divisor graph $\Gamma[\mathbb{Z}_{p^3}]$ is $p-1$ connected.
 \end{proof}
\end{theorem}
With similar arguments, we prove the more general case in the following theorem.
\begin{theorem}
The vertex connectivity of $\Gamma[\mathbb{Z}_{p^n}]$ is $p-1$ $\forall$ $n\geq 3$.
 \begin{proof}
We divide the elements(vertices) of $\Gamma[\mathbb{Z}_{p^n}]$ into $n-1$ disjoint sets namely multiples of $p$, multiples of $p^2$... multiples of $p^{n-1}$, given by
\begin{align*}
 A_{1} &=\{k_{1}p\,|\,k_{1}= 1,2,3,....,p^{n-1}-1\, \text{and}\, p\nmid k_{1} \}  \\
 A_{2} &=\{k_{2}p^2\,|\,k_{2}= 1,2,3,....,p^{n-2}-1\, \text{and}\, p\nmid k_{2} \}\\ 
A_{i} &=\{k_{i}p^i\,|\,k_{i}= 1,2,3,....,p^{n-i}-1\, \text{and}\, p\nmid k_{i} \} 
\end{align*}
Among the above, the smallest set is $A_{n-1}$ of order $p-1$.\\
Now if we leave, even one vertex from $A_{i}$ for $i=1,2,......n-1$ and remove remaining, the graph is still connected because every element of $A_{i}$ is adjacent with $A_{n-1}$.\\
Now if we delete all the vertices of $A_{n-1}$, we get a disconnected graph because the elements of $A_{1}$ are adjacent with only the elements of $A_{n-1}$ , these vertices will become isolated.\\
Therefore the minimum number of vertices to be deleted to make the graph disconnected is $|A_{n-1}|=p-1$.\\
Thus vertex connectivity of $\Gamma[\mathbb{Z}_{p^n}]$ is $p-1$ \\
Since, $\kappa(\Gamma[\mathbb{Z}_{p^n}])=p-1$, the zero divisor graph $\Gamma[\mathbb{Z}_{p^n}]$ is  $p-1$ connected.
 \end{proof}
\end{theorem}
Now before going to the vertex connectivity of  $\Gamma[\mathbb{Z}_{n}]$ for any $n\geq 1$. we first work into the vertex connectivity of  $\Gamma[\mathbb{Z}_{m}]$  where $ m=p^{\alpha}q^{\beta} $.
\begin{theorem}
The vertex connectivity of $\Gamma[\mathbb{Z}_{p^{\alpha}q^{\beta} }]$ is $\min\{p-1,q-1\}$.
 \begin{proof}
Here, we divide the elements(vertices) of $\Gamma[\mathbb{Z}_{m}]$   into  disjoint sets namely multiples of $ p^{i}$, multiples of $q^{j}$ and multiples of $p^{i}q^{j}$ given by
\begin{align*}
 A_{p^{i}} &=\{r_{i}p^{i}\,|\,r_{i}= 1,2,3,....,p^{i}-1\, \text{and}\, p\nmid r_{i} \}\\
 A_{q^{j}} &=\{s_{j}q^{j}\,|\,s_{j}= 1,2,3,....,q^{j}-1\, \text{and}\, q\nmid s_{j} \} \\ 
A_{p^{i}q^{j}} &=\{t_{ij}p^iq^{j}\,|\,t_{ij}= 1,2,3,....,p^{i}q^{j}-1\, \text{and}\, p\nmid t_{ij}\, \text{and}\, q\nmid t_{ij} \}.
\end{align*}
The smallest set among the above is either $B=A_{p^{\alpha}q^{\beta-1}}$ or $C=A_{p^{\alpha-1}q^{\beta}}$.\\
Suppose $p{\displaystyle <\,}q$.\\
Here we note that  every element of $A_{p^{i}}$  is adjacent with each and every element of $C$.\\
So, if we delete all the elements of $C$ then the graph becomes disconnected as the vertices of $A_{p}$ becomes isolated.\\
Therefore the minimum number of vertices to be deleted to make the graph disconnected is $p-1$.\\
Similarly, we get that if $q<p$, then the minimum number of vertices to be deleted to make the graph disconnected is $q-1$.\\
Hence the vertex connectivity of $\Gamma[\mathbb{Z}_{m}]$ is  $\min\{p-1, q-1\}$.
 \end{proof}
\end{theorem}
In the next theorem we derive the vertex connectivity of $\Gamma[\mathbb{Z}_{n }]$ for any natural number $n$.
\begin{theorem}\label{t1}
The vertex connectivity of $\Gamma[\mathbb{Z}_{n }]$ where $n =p_{1}^{\alpha_{1}}p_{2}^{\alpha_{2}}.....p_{k}^{\alpha_{k}}$ is $\min\{p_{1}-1, p_{2}-1,...,p_{k}-1\}$.
 \begin{proof}
consider a zero divisor graph  $\Gamma[\mathbb{Z}_{n}]$  where $n =p_{1}^{\alpha_{1}}p_{2}^{\alpha_{2}}.....p_{k}^{\alpha_{k}}$  .\\
Here,we divide the elements(vertices) of $\Gamma[\mathbb{Z}_{n}]$  into  the corresponding  disjoint sets of product of all possible powers of given primes like set of powers of $p_j^i$, set of product of powers of $p_j^ip_r^s$ and so on.\\
Among these sets, we consider the sets of the form\\
$A_{i} =\{ m( p_{1}^{\alpha_{1}}p_{2}^{\alpha_{2}}.....p_{i}^{\alpha_{i}-1}.....p_{k}^{\alpha_{k}})\,|\,m\nmid p_j\,\text{for all}\,j\}$ with  $|A_{i}|=(p_{i}-1)$.\\
Let $ r=min\{|A_{1}|,...,|A_{k}|\}=min\{p_1-1,...,p_k-1\}$.\\
Suppose $ r=|A_{j}| =p_{j}-1$ ,for some $j$.\\
Now consider the set  $ A_{p_j} =\{ t p_{j}\,|\, p_j\nmid t \}$.\\
since the elements of $ A_{p_j}$ are adjacent only with the vertices of $ A_j $, so deletion of all vertices of $A_{j}$ leaves the elements of $A_{p_j} $  isolated .\\
Thus the graph is disconnected if we remove all the elements of $A_j$.\\
Therefore the minimum number of vertices to be deleted to make the graph disconnected is $r=|A_{j}|$.\\
Thus Vertex connectivity of $\Gamma[\mathbb{Z}_{n}]$ is  $r=\min\{p_{1}-1, p_{2}-1,...,p_{k}-1\}$.
 \end{proof}
\end{theorem}
\section{Edge connectivity of the zero divisor graph  $\Gamma[{\mathbb{Z}_n}]$}
In this section we discuss the edge connectivity of the zero divisor graph $\Gamma[{\mathbb{Z}_n}]$ . To start with we consider $n=p^2$.
\begin{theorem}
The Edge connectivity of $\Gamma[\mathbb{Z}_{p^2}]$ is $p-2$.
 \begin{proof}
The vertex set of the zero divisor graph $\Gamma[\mathbb{Z}_{p^2}]$ is 
$A =\{kp\,|\,k= 1,2,3,....,p-1\, \text{and}\, k\nmid p \}$ and so $|A|=(p-1)$.\\
As the graph is complete, every vertex is incident with $(p-2)$ edges.\\
so removing this (p-2) edges with respect to a fixed vertex $v$, the graph becomes disconnected as $v$ becomes isolated.\\
If we remove fewer than $(p-2)$ edges say $(p-3)$, then as the vertex $v$ is incident with every edge, the vertex $v$ is not isolated so that the graph is still connected.\\
Thus the minimum number of edges to be deleted is $p-2$.\\
Hence the edge connectivity of $\Gamma[\mathbb{Z}_{p^2}]$ is $\kappa_e(\Gamma[\mathbb{Z}_{p^2}])= p-2$.\\
Since $\kappa_e(\Gamma[\mathbb{Z}_{p^2}])= p-2$, the zero divisor graph $\Gamma[\mathbb{Z}_{p^2}]$  is $(p-2)-$ edge connected.
 \end{proof}
\end{theorem}
\begin{theorem}
The Edge connectivity of $\Gamma[\mathbb{Z}_{p^n}]$ is $p-1$ $\forall$ $n\geq 3$. 
 \begin{proof}
We group the elements(vertices) of $\Gamma[\mathbb{Z}_{p^n}]$ into (n-1) disjoint sets namely multiples of p, multiples of $p^2$ and so on multiples of $p^{n-1}$  given by
\begin{align*}
A_{p^{i}} &=\{r_{i}p^{i}\,|\,r_{i}= 1,2,3,....,p^{i}-1\, \text{and}\, p\nmid r_{i} \}\\
A_{q^{j}} &=\{s_{j}q^{j}\,|\,s_{j}= 1,2,3,....,q^{j}-1\, \text{and}\, q\nmid s_{j} \} \\ 
A_{p^{i}q^{j}} &=\{t_{ij}p^iq^{j}\,|\,t_{ij}= 1,2,3,....,p^{i}q^{j}-1\, \text{and}\, p\nmid t_{ij}\, \text{and}\, q\nmid t_{ij} \}.
\end{align*}
with cardinality  $|A_{i}|=(p^{n-i}-1)$ ,for $i=1,2,......n-1$.\\
Clearly, the smallest set among the above is $A_{n-1}$ of length $p-1$.\\
Here $A_{1}$ is the only set in which there exists a vertex (in fact every vertex) which is incident with only $(p-1)$ edges  because any vertex in $A_{i}$ is adjacent with the elements of the sets $A_{n-i}$, $ A_{n-i+1}$,.....$ A_{n-1}$  and so with $p^{i}+p^{i-1}+....+p-i$ edges.\\
Here our idea is to identify a set whose elements are having edges with only $ A_{n-1}$, and the set $A_1$ will do the work for us.\\
Thus if we remove all the edges incident with a vertex  $v$ in $ A_{1}$, which are $(p-1)$ in number as the vertices of $ A_{1}$  are adjacent with only the elements of the set $ A_{n-1}$, giving rise to a disconnected graph.\\  
If we remove fewer than (p-1) edges say (p-2) then  the graph is still connected.\\
Therefore the minimum number of edges to be removed is $p-1$.\\
Thus the edge connectivity of $\Gamma[\mathbb{Z}_{p^n}]$ is $\kappa_{e}(\Gamma[\mathbb{Z}_{p^n}])=p-1$.\\
Since $\kappa_{e}(\Gamma[\mathbb{Z}_{p^n}])=p-1$,  the zero divisor graph $\Gamma[\mathbb{Z}_{p^n}]$ is $(p-1)-$ edge  connected.
 \end{proof}
\end{theorem}
Now we move to the general case: the Edge connectivity of $\Gamma[\mathbb{Z}_{n }]$ for $n =p_{1}^{\alpha_{1}}p_{2}^{\alpha_{2}}.....p_{k}^{\alpha_{k}}$ 
\begin{theorem}
The Edge connectivity of $\Gamma[\mathbb{Z}_{n }]$ for $n =p_{1}^{\alpha_{1}}p_{2}^{\alpha_{2}}.....p_{k}^{\alpha_{k}}$ is r .\\
That is r-edge connected where $ r=min\{p_{i}-1\}$  for $i=1,2,......k$.\\ 
 \begin{proof}
We group the vertex set and choose the suitable set $ A_{j}$ as in Theorem $\ref{t1}$.
Since the vertices of $ A_{p_j}$ are  adjacent only with the vertices of $A_j$, the deletion of  all edges incident with the  any vertex of  $ A_{p_j}$ leaves that vertex isolated and hence the  graph becomes disconnected.
Also in all other possible combinations of $A_{p_s}$  and $A_{s}$ the graph is still connected.\\
Therefore, the minimum number of edges to be deleted to make the graph disconnected is $r=|A_{j}|$.\\
Thus edge connectivity of $\Gamma[\mathbb{Z}_{n}]$ is  $\kappa_{e}[\Gamma[\mathbb{Z}_{n}]=r=\min\{p_{1}-1, p_{2}-1,...,p_{k}-1\}$.
\end{proof}
\end{theorem}
\begin{remark}
 By observing all the results from the sections 3 and 4, we conclude that vertex connectivity is same as the edge connectivity of the zero divisor graph $\Gamma[\mathbb{Z}_{n}]$. i.e. $\kappa(\Gamma[\mathbb{Z}_{n}])= \kappa_{e}(\Gamma[\mathbb{Z}_{n}])$.
\end{remark}
\begin{theorem}
The minimum degree of the zero divisor graph $\Gamma[\mathbb{Z}_{p^n}]$ is equal its edge connectivity. Therefore $\kappa(\Gamma[\mathbb{Z}_{n }])= \kappa_{e}(\Gamma[\mathbb{Z}_{n }]) = \delta(\Gamma[\mathbb{Z}_{n }])$. 
\begin{proof}
Let the edge connectivity of $\Gamma[\mathbb{Z}_{p^n}]$ is $r$.\\
Then the minimum number of edges that are incident to any arbitrary vertex in the graph is $r$.\\
Also there exists a set ( as proved in the previous results) in which degree of every vertex  is $r$ and hence the minimum degree of $\Gamma[\mathbb{Z}_{n}]$ is $\delta(\Gamma(\mathbb{Z}_{n})=r$.\\
Hence $\kappa(\Gamma[\mathbb{Z}_{n }])= \kappa_{e}(\Gamma[\mathbb{Z}_{n }]) = \delta(\Gamma[\mathbb{Z}_{n }])$. 
\end{proof}
\end{theorem}

\end{document}